# SMALL GENERATORS OF NUMBER FIELDS

WOLFGANG M. RUPPERT

1. INTRODUCTION

Let $K$ be a number field of degree $n$ over $\mathbf{Q}$ and $\alpha \in K$ a generator, i.e. $K = \mathbf{Q}(\alpha)$, with minimal polynomial $f = a_n x^n + a_{n-1} x^{n-1} + \cdots + a_1 x + a_0$ and $a_i \in \mathbf{Z}$, $\gcd(a_0, \ldots, a_n) = 1$. We call

$$H(\alpha) = H(f) = \max(|a_0|, \ldots, |a_n|)$$

the height of $\alpha$. Our question is now: How large are small generators $\alpha$ of $K$ (where our measure is the height $H(\alpha)$)?

A bound from below is given in the following

**Proposition 1.** *For every $n \in \mathbf{N}$ there is a real number $c_n > 0$ such that if $\alpha$ generates a number field $K$ of degree $n$ and discriminant $D_K$ then*

$$H(\alpha) \geq c_n |D_K|^{\frac{1}{2n-2}}.$$

*One can take $c_n = \frac{1}{n\sqrt{n}}$.*

*Proof:* Let $f = a_n x^n + \cdots + a_0$ be the minimal polynomial of $\alpha$. According to [Cohen, p.216, exercise 15 (H. W. Lenstra)] the $\mathbf{Z}$-module

$$\begin{aligned} R &= \mathbf{Z} + (a_n \alpha)\mathbf{Z} + (a_n \alpha^2 + a_{n-1}\alpha)\mathbf{Z} + (a_n \alpha^3 + a_{n-1}\alpha^2 + a_{n-2}\alpha)\mathbf{Z} \\ &\quad + \cdots + (a_n \alpha^{n-2} + \cdots + a_3\alpha)\mathbf{Z} + (a_n \alpha^{n-1} + \cdots + a_2\alpha)\mathbf{Z} \end{aligned}$$

is an order in $K = \mathbf{Q}(\alpha)$ with discriminant $D(f)$, the discriminant of the polynomial $f$. So there is an $m \in \mathbf{N}$ with $D(f) = m^2 D_K$. On the other hand $D(f)$ is homogeneous in $a_0, \ldots, a_n$ of degree $2n - 2$ and there is a constant $e_n > 0$ such that

$$|D(f)| \leq e_n H(f)^{2n-2}.$$

This gives $|D_K| \leq e_n H(\alpha)^{2n-2}$ and therefore

$$H(\alpha) \geq e_n^{-\frac{1}{2n-2}} |D_K|^{\frac{1}{2n-2}}.$$

Using the determinant representation of the discriminant one sees that one can take $e_n = (n\sqrt{n})^{2n-2}$, which gives the desired explicit inequality. ∎

A natural question is now:

**Question 1.** *Is there a number $d_n$ such that every number field $K$ of degree $n$ over $\mathbf{Q}$ has a generator $\alpha$ with*

$$H(\alpha) \leq d_n |D_K|^{\frac{1}{2n-2}}?$$

The answer will be 'yes' for quadratic fields. We do not know what's going on in general. We conclude with two examples.

## 2. IMAGINARY QUADRATIC FIELDS

If $K = \mathbf{Q}(\alpha)$ is imaginary quadratic with discriminant $D$ and $f = ax^2 + bx + c$ the minimal polynomial of $\alpha$ then there is an $m \in \mathbf{N}$ with $b^2 - 4ac = D(f) = m^2 \cdot D$ and therefore $4H(\alpha)^2 \geq 4ac = b^2 - m^2 D = b^2 + m^2|D| \geq |D|$, so

$$H(\alpha) \geq \frac{1}{2}\sqrt{|D|}.$$

Is is easy to write down examples which show that the inequality is best possible:

**Example:** Assume that $m \in \mathbf{N}$ and $4m^2 - 1 = (2m-1)(2m+1)$ is square-free. Then $D = 1 - 4m^2$ is the discriminant of $K = \mathbf{Q}(\sqrt{D})$ and is generated by an element $\alpha$ with minimal polynomial $f = mx^2 + x + m$. So we get

$$\frac{H(\alpha)}{\sqrt{|D|}} = \frac{m}{\sqrt{4m^2 - 1}}$$

which tends to $\frac{1}{2}$ as $m$ goes to $\infty$.

To get an estimate in the other direction we look first at small integral elements. The result is given in the following lemma:

**Lemma 1.** *Let $d$ be a square-free integer $< 0$. The smallest integral generators of $\mathbf{Q}(\sqrt{d})$ have minimal polynomial*

$$x^2 - d \quad for \quad d \equiv 2, 3 \bmod 4$$
$$x^2 \pm x + \frac{1-d}{4} \quad for \quad d \equiv 1 \bmod 4.$$

The lemma is easy to prove. It shows that the height of integral generators is always $\geq \frac{1}{4}|D|$ and is far away from what we are looking for.

For a quadratic number field with discriminant $D$ let $H_{min}(D)$ be the height of a generator of minimal height. The following table gives the discriminant $D$ where $\frac{H_{min}(D)}{\sqrt{|D|}}$ is maximal in the specified range. $ax^2 + bx + c$ is the minimal polynomial of a minimal generator of $\mathbf{Q}(\sqrt{D})$.

|  | $D$ | $(a,b,c)$ | $\frac{H_{min}(D)}{\sqrt{|D|}}$ |
|---|---|---|---|
| $0 \geq D \geq -10000$ | -163 | (1,1,41) | 3.2114 |
| $-10000 \geq D \geq -20000$ | -17467 | (47,39,101) | 0.7642 |
| $-20000 \geq D \geq -30000$ | -21379 | (55,29,101) | 0.6908 |
| $-30000 \geq D \geq -40000$ | -36523 | (73,59,137) | 0.7169 |
| $-40000 \geq D \geq -50000$ | -47947 | (83,39,149) | 0.6805 |
| $-50000 \geq D \geq -60000$ | -50395 | (89,35,145) | 0.6459 |
| $-60000 \geq D \geq -70000$ | -68707 | (127,127,167) | 0.6371 |
| $-70000 \geq D \geq -80000$ | -73747 | (109,41,173) | 0.6372 |
| $-80000 \geq D \geq -90000$ | -81859 | (121,93,187) | 0.6536 |
| $-90000 \geq D \geq -100000$ | -91795 | (127,91,197) | 0.6502 |

The table and some further hints suggest the following

**Conjecture 1.** *If $K$ is an imaginary quadratic field with discriminant $D$ then there is $\alpha$ with $K = \mathbf{Q}(\alpha)$ and*

$$H(\alpha) \leq 3.22\sqrt{|D|}.$$

We have the following asymptotic result:

**Theorem 1.**
$$\lim_{D \to -\infty} \frac{H_{min}(D)}{\sqrt{|D|}} = \frac{1}{2}$$
*where $D$ runs through the discriminants of imaginary quadratic fields.*

The proof depends heavily on a result of Duke [Duke] which we shortly recall[1]:
Let
$$F = \{z \in \mathbf{C} : \mathrm{Im} z > 0, -\frac{1}{2} < \mathrm{Re} z \leq \frac{1}{2}, |z| \geq 1 \text{ and } |z| > 1 \text{ if } \mathrm{Re} z < 0\}$$
be the standard fundamental domain of $\mathrm{PSL}_2(\mathbf{Z})$ in the upper half plane $\{z : \mathrm{Im} z > 0\}$. Using the measure $d\mu = \frac{3}{\pi} \frac{dxdy}{y^2}$ gives $\mu(F) = 1$. Define for a discriminant $D$
$$\Lambda_D = \{z = \frac{b + \sqrt{D}}{2a} : b^2 - 4ac = D, z \in F, a, b, c \in \mathbf{Z}\}.$$
If $\Omega \subseteq F$ is convex (in the non-Euclidean sense) with a piece-wise smooth boundary then by [Duke, Theorem 1]
$$\lim_{D \to -\infty} \frac{\#\Lambda_D \cap \Omega}{\#\Lambda_D} = \mu(\Omega)$$
where $D$ runs through the discriminants of imaginary quadratic fields. Now we are ready to prove the theorem:

*Proof of Theorem 1:* Let $0 < \epsilon \leq 1$ be given and choose a convex set $\Omega$ with piece-wise smooth boundary, $\mu(\Omega) > 0$ and
$$\Omega \subseteq \{z \in F : 0 \leq \mathrm{Re} z \leq \frac{1}{2}\epsilon, 1 \leq \mathrm{Im} z \leq 1 + \frac{1}{2}\epsilon\}.$$
By the just mentioned result we find $D_\epsilon$ such that $\#\Lambda_D \cap \Omega \geq 1$ whenever $D < D_\epsilon$. Take $D < D_\epsilon$ then we get an $\alpha \in \Lambda_D \cap \Omega$, i.e. $\alpha = \frac{b+\sqrt{D}}{2a}$ with $a, b, c \in \mathbf{Z}$ and $D = b^2 - 4ac$ such that
$$a > 0, \quad 0 \leq \frac{b}{2a} \leq \frac{1}{2}\epsilon, \quad 1 \leq \frac{\sqrt{|D|}}{2a} \leq 1 + \frac{1}{2}\epsilon.$$
We see at once
$$|a| = a \leq \frac{1}{2}\sqrt{|D|} \quad \text{and} \quad |b| = b \leq a\epsilon \leq a \leq \frac{1}{2}\sqrt{|D|}.$$
Finally
$$\begin{aligned}
|c| &= c = \frac{b^2 + |D|}{4a} \leq \frac{a^2\epsilon^2 + |D|}{4a} = \\
&= a\frac{\epsilon^2}{4} + \frac{1}{2}\sqrt{|D|}\frac{\sqrt{|D|}}{2a} \leq \frac{1}{2}\sqrt{|D|}\frac{\epsilon^2}{4} + \frac{1}{2}\sqrt{|D|}(1 + \frac{1}{2}\epsilon) = \\
&= \frac{1}{2}\sqrt{|D|}(1 + \frac{\epsilon}{2} + \frac{\epsilon^2}{4}) \leq \frac{1}{2}\sqrt{|D|}(1 + \epsilon).
\end{aligned}$$
As $\alpha$ has minimal polynomial $ax^2 - bx + c$ and generates $\mathbf{Q}(\sqrt{D})$ we get (with the trivial estimate)
$$\frac{1}{2} \leq \frac{H_{min}(D)}{\sqrt{|D|}} \leq \frac{H(\alpha)}{\sqrt{|D|}} \leq \frac{1}{2}(1 + \epsilon)$$

---

[1] I would like to thank H.W.Lenstra who gave me the hint to Duke's paper

which proves our claim. ∎

It would be nice to have an effective version of Duke's theorem in order to prove a statement like Conjecture 1.

## 3. Real quadratic fields

Let $K$ be a real quadratic field with discriminant $D$ and $\alpha$ a generator of $K$ with minimal polynomial $f = ax^2 + bx + c$. Then there is an $m \in \mathbf{N}$ with $b^2 - 4ac = m^2 D$ so that we obtain $0 < D \leq m^2 D = b^2 - 4ac \leq b^2 + 4|a||c| \leq 5H(\alpha)^2$ and therefore

$$H(\alpha) \geq \frac{1}{\sqrt{5}}\sqrt{D} = 0.4472\ldots\sqrt{D}.$$

The following example shows that the estimate is best possible:

**Example:** Let $m \in \mathbf{N}$ such that $5m^2 - 2m + 1$ is square-free. Then an element $\alpha$ with minimal polynomial $f = mx^2 + (m-1)x - m$ generates a real quadratic number field $K$ with discriminant $D = 5m^2 - 2m + 1$ and

$$\frac{H(\alpha)}{\sqrt{D}} = \frac{m}{\sqrt{5m^2 - 2m + 1}}$$

tends to $\frac{1}{\sqrt{5}}$ as $m$ goes to $\infty$.

To get an estimate in the other direction the situation is much easier than in the imaginary quadratic case as we find small integral generators:

**Proposition 2.** *If $K$ is a real quadratic field with discriminant $D$ there is an (integral) $\alpha \in K$ with $K = \mathbf{Q}(\alpha)$ and*

$$H(\alpha) < \sqrt{D}.$$

*Proof:* Let $m = [\sqrt{D}]$ and choose $a = m$ or $a = m - 1$ such that $a^2 \equiv D \mod 4$. Take $b$ with $a^2 - 4b = D$. It is clear that $b < 0$. The assumption $|b| \geq m$ would imply

$$D = a^2 + 4|b| \geq (m-1)^2 + 4m = m^2 + 2m + 1 = (m+1)^2,$$

so $\sqrt{D} \geq m+1$, which contradicts the definition of $m$. Therefore an element $\alpha$ with minimal polynomial $f = x^2 + ax + b$ has height $< \sqrt{D}$ and proves the proposition. ∎

Proposition 2 answers Question 1 for real quadratic fields in an effective way. For the rest of this section we want to study the asymptotic behavior of $\frac{H_{min}(D)}{\sqrt{D}}$. We know already

$$\frac{1}{\sqrt{5}} \leq \frac{H_{min}(D)}{\sqrt{D}} < 1 \quad \text{and} \quad \liminf_{D \to \infty} \frac{H_{min}(D)}{\sqrt{D}} = \frac{1}{\sqrt{5}}.$$

The following table lists discriminants $D$ where $\frac{H_{min}(D)}{\sqrt{D}}$ is maximal in the specified range. $ax^2 + bx + c$ is the minimal polynomial of a minimal generator.

|  | $D$ | $(a,b,c)$ | $\frac{H_{min}(D)}{\sqrt{D}}$ |
|---|---|---|---|
| $0 \leq D \leq 1000$ | 293 | (1,15,-17) | 0.9932 |
| $1000 \leq D \leq 2000$ | 1592 | (2,36,-37) | 0.9273 |
| $2000 \leq D \leq 3000$ | 2540 | (10,30,-41) | 0.8135 |
| $3000 \leq D \leq 4000$ | 3053 | (7,43,-43) | 0.7782 |
| $4000 \leq D \leq 5000$ | 4973 | (17,37,-53) | 0.7516 |
| $5000 \leq D \leq 6000$ | 5885 | (13,55,-55) | 0.7170 |
| $6000 \leq D \leq 7000$ | 6341 | (17,51,-55) | 0.6907 |
| $7000 \leq D \leq 8000$ | 7229 | (17,53,-65) | 0.7645 |
| $8000 \leq D \leq 9000$ | 8197 | (23,49,-63) | 0.6959 |
| $9000 \leq D \leq 10000$ | 9037 | (37,3,-61) | 0.6417 |

We will study $\limsup_{D \to \infty} \frac{H_{min}(D)}{\sqrt{D}}$ in two ways each of which depends on certain conjectures.

3.1. **Assuming the existence of primes with certain properties.** Let $M_\epsilon$ be the set of all discriminants $D$ of real quadratic fields $\mathbf{Q}(\sqrt{D})$ such that there is an odd prime $p$ with

$$(\frac{D}{p}) = 1 \quad \text{and} \quad \frac{1}{2}\sqrt{D} \leq p \leq (\frac{1}{2} + \epsilon)\sqrt{D}.$$

Then we get the

**Lemma 2.** *If $D \in M_\epsilon$ and $0 < \epsilon \leq \frac{1}{2}$ the real quadratic field $\mathbf{Q}(\sqrt{D})$ has a generator $\alpha$ with*

$$H(\alpha) \leq (\frac{1}{2} + \epsilon)\sqrt{D}.$$

*Proof:* By definition there is an odd prime $p$ with

$$(\frac{D}{p}) = 1 \quad \text{and} \quad \frac{1}{2}\sqrt{D} \leq p \leq (\frac{1}{2} + \epsilon)\sqrt{D} \leq \sqrt{D}.$$

Choose $b \in \mathbf{Z}$ with $0 \leq b \leq p$ and $b^2 \equiv D \bmod p$. The number $p-b$ satisfies the same conditions so that we can assume $b^2 \equiv D \bmod 2p$ which implies $b^2 \equiv D \bmod 4p$ as $D \equiv 0, 1 \bmod 4$. Define $c \in \mathbf{Z}$ as $c = \frac{b^2-D}{4p}$. As $b^2 < D$ we have $c < 0$. It follows

$$|c| = -c = \frac{D-b^2}{4p} \leq \frac{D}{4p} \leq \frac{1}{2}\sqrt{D}.$$

The element $\alpha$ with minimal polynomial $f = px^2 + bx + c$ generates $\mathbf{Q}(\sqrt{D})$ and satisfies

$$H(\alpha) \leq (\frac{1}{2} + \epsilon)\sqrt{D}$$

which proves the lemma. ∎

If $\epsilon > 0$ is fixed and $D$ is large there are many primes $p$ with $\frac{1}{2}\sqrt{D} \leq p \leq (\frac{1}{2}+\epsilon)\sqrt{D}$. As $(\frac{D}{p})$ is a quadratic character the following conjecture seems plausible:

**Conjecture 2.** *For every $\epsilon > 0$ there is a $D_\epsilon$ such that $D > D_\epsilon$ implies $D \in M_\epsilon$ if $D$ is a real quadratic discriminant.*

After some computations we conjecture e.g. $D > 981913 \Rightarrow D \in M_{0.1}$.

I do not know how to attack conjecture 2. A natural question seems to be:

**Question 2.** *Can conjecture 2 be deduced from ERH (extended Riemann hypothesis)?*

An immediate consequence of Lemma 2 is

**Corollary 1.** *If conjecture 2 holds then*
$$\limsup_{D\to\infty} \frac{H_{min}(D)}{\sqrt{D}} \leq \frac{1}{2}$$
*where $D$ runs through the discriminants of real quadratic fields.*

### 3.2. Assuming the existence of reduced elements with certain properties.

Let $D$ be the discriminant of a real quadratic field and define the set of reduced elements by
$$\Lambda_D = \{\alpha = \frac{b+\sqrt{D}}{2a} : b^2 - 4ac = D, \alpha > 1, -1 < \alpha' < 0, a, b, c \in \mathbf{Z}\}$$
(where $\alpha'$ is the conjugate of $\alpha$ and $\sqrt{D} > 0$). The map $\alpha \mapsto \frac{1}{\alpha - [\alpha]}$ induces a decomposition of $\Lambda_D$ in $h(D)$ cycles where $h(D)$ is the class number of $\mathbf{Q}(\sqrt{D})$ [Cohen, p.260]. It is also known that $\lim_{D\to\infty} \frac{\log \#\Lambda_D}{\log \sqrt{D}} = 1$ [Lachaud] so the elements of $\Lambda_D$ play a similar role as the elements of $\Lambda_D$ in the imaginary quadratic case.

The smallest generators of a field need not be reduced as the following example shows:

**Example:** The real quadratic field $K = \mathbf{Q}(\sqrt{635})$ with discriminant $D = 2540$ has as a smallest generator $\alpha = \frac{15+\sqrt{635}}{10}$ with minimal polynomial $10x^2 - 30x - 41$ and height $H(\alpha) = 41$. Among the reduced elements $\beta = 25 + \sqrt{635}$ (with minimal polynomial $x^2 - 50x - 10$) has the smallest height $H(\beta) = 50$.

Nevertheless we have the following lemma:

**Lemma 3.** *If $\alpha$ generates $K$ with discriminant $D$ and $H(\alpha) \leq 0.48\sqrt{D}$ then one of the elements $\alpha, \alpha', -\alpha, -\alpha'$ is reduced.*

*Proof:* Let $ax^2 - bx + c$ be the minimal polynomial of $\alpha$ with $a > 0$. Changing from $\alpha$ to $-\alpha$ we can assume $b \geq 0$. Then we there is $f \in \mathbf{N}$ with $b^2 - 4ac = Df^2$. As $a, b, |c| \leq 0.48\sqrt{D}$ we get $Df^2 \leq 5 \cdot 0.48^2 D$ and therefore $f = 1$. So we have without restriction
$$\alpha = \frac{b+\sqrt{D}}{2a} \quad \text{and} \quad \alpha' = \frac{b-\sqrt{D}}{2a}.$$
As $D = b^2 - 4ac \leq b^2 + 4 \cdot 0.48^2 D$ we get $b \geq 0.28\sqrt{D}$ and $D \leq 0.48^2 D + 4 \cdot 0.48\sqrt{D} \cdot a$ gives $a \geq 0.40\sqrt{D}$. Therefore
$$\alpha = \frac{b+\sqrt{D}}{2a} \geq \frac{0.28\sqrt{D}+\sqrt{D}}{2 \cdot 0.48\sqrt{D}} \geq 1.33 > 1$$
and
$$0 > \alpha' = \frac{b-\sqrt{D}}{2a} \geq \frac{0.28\sqrt{D}-\sqrt{D}}{2a} = -\frac{0.72\sqrt{D}}{2a} \geq -\frac{0.72\sqrt{D}}{2 \cdot 0.40\sqrt{D}} > -1$$
which shows that $\alpha$ is reduced. ∎

**Lemma 4.** *Let $K$ be a real quadratic field with discriminant $D$. If $\alpha \in \Lambda_D$ then*

$$\frac{H(\alpha)}{\sqrt{D}} = \max(\frac{1}{\alpha - \alpha'}, \frac{\alpha + \alpha'}{\alpha - \alpha'}, \frac{\alpha(-\alpha')}{\alpha - \alpha'}).$$

*Proof:* If the minimal polynomial of $\alpha$ is $ax^2 - bx + c$ then

$$|a| = a = \sqrt{D} \cdot \frac{1}{\alpha - \alpha'},$$

$$|b| = b = \sqrt{D} \cdot \frac{\alpha + \alpha'}{\alpha - \alpha'},$$

$$|c| = -c = \sqrt{D} \cdot \frac{\alpha(-\alpha')}{\alpha - \alpha'},$$

which gives the result. ∎

Let $G = \{(x, y) \in \mathbf{R}^2 : x > 1, -1 < y < 0\}$ and $\tilde{\Lambda}_D = \{(\alpha, \alpha') : \alpha \in \Lambda_D\}$, then $\tilde{\Lambda}_D \subseteq G$. Define for $0 < h < 1$ the set

$$G_h = \{(x, y) \in G : \frac{1}{x - y} \leq h, \frac{x + y}{x - y} \leq h, \frac{x(-y)}{x - y} \leq h\},$$

which can be written as

$$G_h = \{(x, y) \in G : y \leq x - \frac{1}{h}, y \leq -\frac{1 - h}{1 + h}x, y \geq -\frac{hx}{x - h}\}.$$

If $(\alpha, \alpha') \in \tilde{\Lambda}_D \cap G_h$ then the lemma gives $H(\alpha) \leq h\sqrt{D}$. It is not difficult to see that $G_h = \emptyset$ for $h < \frac{1}{\sqrt{5}}$, $G_{\frac{1}{\sqrt{5}}} = \{(\frac{1+\sqrt{5}}{2}, \frac{1-\sqrt{5}}{2})\}$ and $G_h$ is the closure of an open non empty set for $h > \frac{1}{\sqrt{5}}$.

Now we formulate a conjecture for reduced real quadratic numbers:

**Conjecture 3.** *If $U$ is an open non empty subset of $G$ then there is a $c_U$ such that*

$$U \cap \tilde{\Lambda}_D \neq \emptyset \text{ for all } D > c_U.$$

It would be interesting to know if Conjecture 3 can be deduced from Duke's results for real quadratic fields [Duke].

**Corollary 2.** *Assuming Conjecture 3 we get*

$$\lim_{D \to \infty} \frac{H_{min}(D)}{\sqrt{D}} = \frac{1}{\sqrt{5}}$$

*where $D$ runs through the discriminants of real quadratic fields.*

*Proof:* Let $\epsilon > 0$ be sufficiently small. Then there is a $d_\epsilon$ such that $G_{\frac{1}{\sqrt{5}} + \epsilon} \cap \tilde{\Lambda}_D \neq \emptyset$ for all discriminant $D > d_\epsilon$. For $D > d_\epsilon$ take $\alpha \in \Lambda_D$ with $(\alpha, \alpha') \in \tilde{\Lambda}_D \cap G_{\frac{1}{\sqrt{5}} + \epsilon}$. Then

$$\frac{1}{\sqrt{5}} \leq \frac{H(\alpha)}{\sqrt{D}} \leq \frac{1}{\sqrt{5}} + \epsilon,$$

which implies at once our statement. ∎

We conclude this section with numerical examples. Let $H_{min,red}(D)$ be the minimal height of all elements in $\Lambda_D$. The following table gives the discriminant $D$ where $\frac{H_{min,red}(D)}{\sqrt{D}}$ is maximal in the specified range. $ax^2 - bx + c$ is the minimal polynomial of a corresponding element of $\Lambda_D$. Finally *average* gives the average value of all $\frac{H_{min,red}(D')}{\sqrt{D'}}$ in the given range.

|  | $D$ | $(a,b,c)$ | $\frac{H_{min,red}(D)}{\sqrt{D}}$ | average |
|---|---|---|---|---|
| $1 \leq D \leq 10000$ | 908 | (1,30,-2) | 0.9956 | 0.5238 |
| $10000 \leq D \leq 20000$ | 14693 | (19,109,-37) | 0.8992 | 0.4976 |
| $20000 \leq D \leq 30000$ | 24173 | (23,115,-119) | 0.7654 | 0.4904 |
| $30000 \leq D \leq 40000$ | 37532 | (38,122,-149) | 0.7691 | 0.4881 |
| $40000 \leq D \leq 50000$ | 49013 | (37,153,-173) | 0.7814 | 0.4847 |
| $50000 \leq D \leq 60000$ | 54053 | (47,153,-163) | 0.7011 | 0.4836 |
| $60000 \leq D \leq 70000$ | 69893 | (97,173,-103) | 0.6544 | 0.4820 |
| $70000 \leq D \leq 80000$ | 79805 | (95,105,-181) | 0.6407 | 0.4814 |
| $80000 \leq D \leq 90000$ | 87533 | (79,159,-197) | 0.6659 | 0.4801 |
| $90000 \leq D \leq 100000$ | 95672 | (106,128,-187) | 0.6046 | 0.4794 |
| $10^5 \leq D \leq 10^5 + 10000$ | 104093 | (83,195,-199) | 0.6168 | 0.4791 |
| $10^6 \leq D \leq 10^6 + 10000$ | 1006232 | (463,194,-523) | 0.5214 | 0.4668 |
| $10^7 \leq D \leq 10^7 + 10000$ | 10000973 | (1423,1029,-1571) | 0.4968 | 0.4589 |

## 4. Concluding Remarks

The following example shows that for every $n$ there are infinitely many number fields of degree $n$ such that an estimate as in question 1 holds.

**Example:** Let $n$ be an integer $\geq 2$ and $p, q$ primes with $p < q < 2p$. Let $\alpha$ be a zero of $f = px^n + q$ and $K = \mathbf{Q}(\alpha)$. Then $K$ has degree $n$ over $\mathbf{Q}$ and the primes $p$ and $q$ are totally ramified in $K$, so $p^{n-1}$ and $q^{n-1}$ divide $D_K$. Therefore we get the estimate

$$H(\alpha) = q < \sqrt{2pq} = \sqrt{2}(p^{n-1}q^{n-1})^{\frac{1}{2n-2}} \leq \sqrt{2}|D_K|^{\frac{1}{2n-2}}.$$

In the next example 'small' (integral) generators are constructed with Minkowski's theorem.

**Example:** Let $K$ be a totally real number field of prime degree $n$, $\alpha_1, \ldots, \alpha_n$ an integral basis of $K$ and $\sigma_1, \ldots, \sigma_n$ the different embeddings $K \hookrightarrow \mathbf{R}$. As $|\det(\sigma_i\alpha_j)_{i,j}| = \sqrt{|D_K|}$ there is $(x_1, \ldots, x_n) \in \mathbf{Z}^n \setminus \{0\}$ by Minkowski's linear forms theorem [Hua, p.540] such that $\alpha = x_1\alpha_1 + \cdots + x_n\alpha_n$ satisfies

$$|\sigma_1\alpha| < 1 \quad \text{and} \quad |\sigma_2\alpha|, \ldots, |\sigma_n\alpha| \leq |D_K|^{\frac{1}{2n-2}}.$$

The condition $|\sigma_1\alpha| < 1$ implies $\alpha \notin \mathbf{Z}$, so $K = \mathbf{Q}(\alpha)$. Let $f = x^n + a_1x^{n-1} + \cdots + a_n$ be the minimal polynomial of $\alpha$. Then for $1 \leq j \leq n-1$ we get $|a_j| \leq \binom{n}{j}|D_K|^{\frac{j}{2n-2}} \leq 2^n\sqrt{|D_K|}$ and $|a_n| \leq \sqrt{|D_K|}$, therefore

$$H(\alpha) \leq 2^n\sqrt{|D_K|}.$$

In case $n \geq 3$ this is far away from what we would like to have.

MATHEMATISCHES INSTITUT, UNIVERSITÄT ERLANGEN-NÜRNBERG, BISMARCKSTRASSE $1\frac{1}{2}$, D-91054 ERLANGEN, GERMANY

*E-mail address*: `ruppert@mi.uni-erlangen.de`